\documentclass{article}
\usepackage{amsfonts,amsmath, amsthm}
\usepackage{amssymb}
\usepackage[utf8]{inputenc}
\usepackage{graphicx}
\usepackage{mathtext}
\usepackage[T2A]{fontenc}
\usepackage{multirow}
\usepackage{setspace}
\usepackage{dsfont}

\hoffset= - 1.1in

\theoremstyle{definition}
\newtheorem{prob}{}
\numberwithin{equation}{subsection}

\newtheorem{theorem}{Theorem}[section]

\newtheorem{lemma}[theorem]{Lemma}

\newenvironment{my_enumerate}{
\vspace{-7pt}
\begin{enumerate}
  \setlength{\itemsep}{1pt}
  \setlength{\parskip}{0pt}
  \setlength{\parsep}{0pt}}{\end{enumerate}

\vspace{-10pt}
}

\newcommand{\eq}[1]{\begin{equation} #1 \end{equation}}

\newcommand{\br}[1]{\left( #1 \right)}

\renewcommand{\l}{\lambda}

\newcommand{\eqn}[1]{\begin{equation*} #1 \end{equation*}}

\begin{document}

\title{
\begin{flushright}
 \mbox{\normalsize ITEP/TH-48/08}
\end{flushright}
\vskip 20pt
On the three-term curves
}

\author{George B. Shabat\footnote{E-mail address: george.shabat@gmail.com, ITEP and RSUH}\\Alexei Sleptsov\footnote{E-mail address: sleptsov@itep.ru, ITEP}} 
\date{}
\maketitle

\begin{abstract}
We classify  projective plane nonsingular curves admitting a 3-term presentation; they exist in any degree, generally constitute 5 birational families and are defined over rational numbers. The Belyi functions on all these curves are presented.
\end{abstract}

\tableofcontents

\section{Introduction}
The general informal problem we are considering in the present paper is:
\begin{center}
WHAT MAKES AN ALGEBRAIC CURVE REMARKABLE?
\end{center}
Some of the curves are so famous that are attributed to great mathematicians;
folium of Descartes, limacon of Pascal, Fermat curves and Klein quartic provide
the obvious examples.

In the cases of such "personal"\  curves a short name defines the curve completely; we
generalize this phenomenon and take the point of view according to which
for a curve \it being remarkable \rm means \it the possibility of defining it
by a short text. \rm \\

This point of view leads to the following question:
\begin{center}
WHAT AMOUNT OF INFORMATION IS NEEDED\\ TO DEFINE AN ALGEBRAIC CURVE?
\end{center}
Perhaps, this question makes sense only over countable fields; actually,
the precise mathematical results of the present paper concern the
curves over $\overline{\mathbb{Q}}$.\\

There are several ways of assigning a mathematical meaning to the
above question. The \emph{Kolmogorov complexity} approach (see, e.g.,\cite{ZvLev}),
as far as the authors know, has not been used so far. The \it heights \rm
on the modular spaces $\mathcal{M}_g(\overline{\mathbb{Q}})$ provide another tool
for measuring the above amount of information; see, e.g., \cite{Liu} and \cite{Lang}.
The \it J.E. Cremona list \rm \cite{Crem} gives
an example of using the modern informational technologies for ordering  great amount
of curves -- namely, the elliptic modular ones -- by some kind of complexity. Finally,
due to the hard \emph{Belyi theorem} (see, e.g., \cite{LandoZv}), one can define the
\emph{Belyi height} of a curve over $\overline{\mathbb{Q}}$ as the smallest degree of a
nonconstant rational function on it with only three critical values; see
\cite{Catalog} for the complete list of curves of Belyi height $\le4$.\\

In the present paper we take the most simple-minded approach to the question, restricting
our attention to \emph{smooth plane projective curves}. It is based on two observations:\\
\\
$\bullet$ all the elliptic curves over $\mathbb{C}$ admit the Weierstrass equation
$$
y^2=4x^3-g_2x-g_3
$$
with four terms, while exactly two of them, with $g_2=0$ and $g_3=0$, are defined by
the three-term equation, and they are doubtlessly
remarkable;\\
\\
$\bullet$ the two remarkable quartics, namely Fermat and Klein ones, are defined by the
three-term equations
$$
x^4+y^4=z^4
$$
and
$$
x^3y+y^3z+z^3x=0.
$$

Thus we are led to the problem of the classification of the \emph{three-term smooth plane
projective curves} up to the birational equivalence; this problem is completely solved
in the present paper, see Theorem \ref{theorem}. It turns out that for any degree $d\ge3$
there exists no more than five
such curves, that we call (by obvious reasons) belonging to the \emph{Fermat} type,
\emph{generalized Klein} type, \emph{block} type, \emph{small Jordan} type and
\emph{big Jordan} type. For $d=3$, as pointed above, they split into two birational classes,
see \ref{cubics}; for $d=4$ the Fermat type coincides with the Big Jordan type; see \ref{quartics}. All the
remaining ones are presumably birationally different.\\

We show in the section \ref{Belyi} that all our curves have low Belyi height. The corresponding
\emph{dessins d'enfants} (see \cite{LandoZv}) will be considered in our future papers, as
well as some obvious generalizations of the problems, considered in the present one.

\section{Three-term curves}

 \normalsize
\subsection{Posing the problem} We take the most
straightforward approach to the above question, considering the
equations of plane projective curves with the
smallest possible number of terms.\\

One checks immediately that the one- and two-termed equations
provide only lines and reducible curves. So we arrive at the main
question of the present paper. Fixing an algebraically closed field
$\Bbbk$ and a \it degree \rm $d$, we pose the following
problem.\\

\it Classify the smooth irreducible curves in $\mathbf{P}^2(\Bbbk)$,
that can be defined in some coordinates $(x:y:z)$ by the equations
of the form
$$
A_1x^{p_{11}}y^{p_{12}}z^{p_{13}}+A_2x^{p_{21}}y^{p_{22}}z^{p_{23}}
+A_3x^{p_{31}}y^{p_{32}}z^{p_{33}} =0,
$$
with $(A_1:A_2:A_3)\in\mathbf{P}^2(\Bbbk)$. \rm\\
\\
Due to the above remarks, we can assume $A_1A_2A_3\ne0$.\\
\\
\subsection{Coefficients}
Introduce \textit{the power matrix}
$$ P:= \left (
 \begin{array}{ccc}
  p_{11}&p_{12}&p_{13}\\
  p_{21}&p_{22}&p_{23}\\
  p_{31}&p_{32}&p_{33}\\
 \end{array}
\right )
$$
 As we are going to see, we can set all of
them to be 1's.\\

\bf Proposition. \it If the power matrix of a three-term  curve
 is non-degenerate, then it  is projectively
equivalent to the one with $A_1=A_2=A_3=1$, i.e. to a curve, defined
by the equation
$$
x^{p_{11}}y^{p_{12}}z^{p_{13}}+x^{p_{21}}y^{p_{22}}z^{p_{23}}
+x^{p_{31}}y^{p_{32}}z^{p_{33}} =0.
$$
\bf Proof. \rm Substituting $x \rightarrow \lambda_1 x$, $y \rightarrow \lambda_2 y$, $z \rightarrow \lambda_3 z$ we want to get
$$
A_1\lambda_1^{p_{11}}\lambda_2^{p_{12}}\lambda_3^{p_{13}}=
A_2\lambda_1^{p_{21}}\lambda_2^{p_{22}}\lambda_3^{p_{23}}=
A_3\lambda_1^{p_{31}}\lambda_2^{p_{32}}\lambda_3^{p_{33}}.
$$
Assume that all these quantities equal 1 and get the system \eqn{
\left\{
\begin{aligned}
& \l_{1}^{p_{11}}\l_2^{p_{12}}\l_3^{p_{13}}=\dfrac{1}{A_1} \ \ \  \ \ \ (1) \\
& \l_{1}^{p_{21}}\l_2^{p_{22}}\l_3^{p_{23}}=\dfrac{1}{A_2}  \ \ \ \ \ \ (2)\\
& \l_{1}^{p_{31}}\l_2^{p_{32}}\l_3^{p_{33}}=\dfrac{1}{A_3} \ \ \ \ \ \ (3). \\
\end{aligned}
\right. }
 Since we suppose that the power matrix $P$ is non-degenerate with
 the \emph{integer} determinant
 $$
 \Delta:=\det P\ne0
 $$
 there is the \emph{adjoint} matrix $Q=(q_{ij})$ satisfying
 $$
 \sum_{i=1}^3q_{ki}p_{ij}=\delta_{ij}\Delta.
 $$
 Then for $k=1,2,3$ we multiply the above equations (1), (2), (3) and obtain
 $$
 \prod_{1\le
 i,j\le3}\lambda_j^{p_{ij}q_{kj}}=
 \prod_{j=1}^3\lambda_j^{\sum_{i=1}^3q_{ki}p_{ij}}=\lambda_k^\Delta=
  \prod_{i=1}^3\frac{1}{A_j^{q_{ki}}}=:B_k,
 $$
 so we take for $\lambda_k$'s any solution of the equation
 $\lambda_k^\Delta=B_k$.
\begin{flushright} \bf QED
\end{flushright}

\subsection{Properties of the power matrices}

\begin{lemma}\label{lemma2} \textit{For the power matrix $p_{ij}$ of irreducible smooth projective curve the following holds:}
\begin{enumerate}
 \item \textit{every column contains zero;}
 \item \textit{every row contains zero;}
 \item \textit{if some row contains the only zero then it contains 1}.
\end{enumerate}
\end{lemma}

\bf Proof.\rm

\

\begin{enumerate}
 \item It obviously follows from the condition of the irreducibility of the curve.

Hence, there is two possibilities: (a) a row contains two zeros; (b) every line contains the only zero. Now we prove that the second and the third statements of the lemma are correct.
 \item
 \begin{enumerate}
   \item A row contains two zeros.

It means that there exists a row which doesn not contain zero. To speak definitely let $p_{11}=p_{12}=0, \ p_{13}=d$. Then the third column contains at least one zero with agreement to the 1). If it contains two zeros then the statement is proved. Because of this let the third column contains the only zero, for expample, $p_{23}=0$ and all other $p_{ij}\neq 0.$ Then the curve takes the form:
\eq{
z^d +  x^{p_{21}} y^{p_{22}} +  x^{p_{31}} y^{p_{32}} z^{p_{33}} = 0.
}
Considering the partial derivatives
\eq{
\left\{
\begin{aligned}
& p_{21}x^{p_{21}-1}y^{p_{22}}+ p_{31}x^{p_{31}-1}y^{p_{32}}z^{p_{33}} = 0 \\
& p_{22}x^{p_{21}}y^{p_{22}-1}+ p_{32}x^{p_{31}}y^{p_{32}-1}z^{p_{33}} = 0 \\
& d\cdot z^{d-1}+p_{33}x^{p_{31}}y^{p_{32}}z^{p_{33}-1} = 0 \\
\end{aligned}
\right.
}
we obtain that $(1:0:0)$ or $(0:1:0)$ is a singular point (some information about singular curves can found in \cite{Shaf}).

 \item Every line contains the only zero

To speak definitely $p_{11}=p_{22}=p_{33}=0$. The curve takes the form
\eq{
y^{p_{12}}z^{p_{13}}+x^{p_{21}}z^{p_{23}}+x^{p_{31}}y^{p_{32}}=0.
}
Let us consider the partial derivatives
\eq{
\left\{
\begin{aligned}
& p_{21}x^{p_{21}-1}z^{p_{23}}+ p_{31}x^{p_{31}-1}y^{p_{32}} = 0 \\
& p_{12}y^{p_{12}-1}z^{p_{13}}+ p_{32}x^{p_{31}}y^{p_{32}-1} = 0 \\
& p_{13}y^{p_{12}}z^{p_{13}-1}+p_{23}x^{p_{12}}z^{p_{23}-1} = 0. \\
\end{aligned}
\right.
}
The points $(1:0:0),(0:1:0),(0:0:1)$ are singular if every row does not contains 1.
\begin{flushright} \bf QED
\end{flushright}

\end{enumerate}
\end{enumerate}

%\bf1.2.1. \it The power matrix contains no zero rows. \rm Indeed,
%the power matrix with a zero row would define a \it two-term \rm
%curve, which, if irreducible, has an affine model $y^m=x^n$ which is
%\it usually \rm singular... (BUT WHAT ABOUT $y=x^n$???)\\
%\\
%\bf1.2.2. \it The power matrix contains no zero columns. \rm Indeed,
%a zero column would mean that one of $x,y,z$ is out of the game, hence
%the curve is a union of lines.\\
%\\
%\bf1.2.3. \it The power matrix contains no nonzero rows. \rm Indeed,
%say, the third nonzero row would mean that $z$ is present in a
%positive degree in all the three terms and hence the curve contains
%a component $z=0$.\\
%\\
%\bf1.2.4. \it The power matrix contains no nonzero columns. \rm
%Indeed, if, say, $a_1b_1c_1\ne0,$ then the only possibility for a curve
%to be irreducible is to be defined by the equation
%$$
%x^{a_1}y^{b_1}z^{c_1}+x^d+y^d=0,
%$$
%but then the point $(0:0:1)$ is singular.
%\\
%\\
\subsection{The main result} We give a complete classification of
 the curves we are considering.\\
 \\

\begin{theorem}\label{theorem}
\textit{Any irreducible projective smooth three-term curve is equivalent to one
 of
 the following five types (where *'s stand for the nonzero elements):}
\end{theorem}
 \begin{center}
\begin{tabular}{|c|c|c|}
\hline
Matrix& Equation & Name \\
\hline
 $\left (
 \begin{array}{ccc}
  *&0&0\\
  0&*&0\\
  0&0&*\\
 \end{array}
\right )$&  $x^d+y^d+z^d=0$ & Fermat (diagonal) type  \\

$\left (
 \begin{array}{ccc}
  *&*&0\\
  0&*&0\\
  0&0&*\\
 \end{array}
\right )$& $xy^{d-1}+y^d+z^d=0$ & Small Jordan type  \\
$\left (
 \begin{array}{ccc}
  *&0&0\\
  0&*&*\\
  0&*&*\\
 \end{array}
\right )$&  $x^d+yz^{d-1}+y^{d-1}z=0$
& block type \\
$ \left (
 \begin{array}{ccc}
  *&*&0\\
  0&*&*\\
  0&0&*\\
 \end{array}
\right )$& $x^d+xy^{d-1}+yz^{d-1}=0$
& Big Jordan type \\
$\left (
 \begin{array}{ccc}
  *&*&0\\
  *&0&*\\
  0&*&*\\
 \end{array}
\right )$&  $xy^{d-1}+yz^{d-1}+zx^{d-1}=0$
  & Klein type \\
 \hline
\end{tabular}
\end{center}
\rm
\bf Proof. \rm Consider our matrices\\
$$
\left (
 \begin{array}{ccc}
  p_{11}&p_{12}&p_{13}\\
  p_{21}&p_{22}&p_{23}\\
  p_{31}&p_{32}&p_{33}\\
 \end{array}
\right )
$$
with $p_{ij}\in\{0,*\}$. There are 512 such
matrices.\\
\\
We consider our matrices up to the action of
$$
\mathbb{S}_3\times\mathbb{S}_3,
$$
permuting rows and columns. Permuting rows corresponds to
rearranging the summands and permuting columns - to renaming
variables.\\

All the 512 matrices split into 32
$\mathbb{S}_3\times\mathbb{S}_3$-orbits.\\

Using lemma \ref{lemma2}, we delete the orbits containing zero rows
and zero columns, and only 16 orbits remain .\\

Using lemma \ref{lemma2}, we delete the orbits containing nonzero
rows and nonzero columns, and only 6 orbits remain.\\

 One of the remaining matrices is impossible. It is
$$
\left (
 \begin{array}{ccc}
  0&*&*\\
  *&0&0\\
  *&0&0\\
 \end{array}
\right )
$$
The corresponding curve is defined by the equation
$$
y^{b_1}z^{c_1}+2x^d=0
$$
and has a singular point $(0:0:1)$.\\

The remaining five matrices constitute the list in the formulation
of the theorem.
\begin{flushright}
\bf QED
\end{flushright}

\subsection{ Non-degeneracy}
\label{nondeg}
 Now it is easy to prove the following\\

\bf Proposition. \it The power matrix of an irreducible smooth three-term
curve is non-degenerate.\\

\bf Proof. \rm Among the five types of matrices, provided by 1.4,
four are of the upper-triangular type, for which the non-degeneracy
is immediate. So only the block type remains, and we have to check
only the non-degeneracy of the $2\times2$-block. It has the form
$$
\left (
 \begin{array}{cc}
  m_1&d-m_2\\
  d-m_1&m_2\\
 \end{array}
\right )
$$
with the determinant $-d^2+d(m_1+m_2)=d(-d+m_1+m_2)$, so the
degeneracy of the power matrix would imply $m_1+m_2=d$, hence the
two equal columns. The curve then will have the form
$$
x^d+2y^{m_1}z^{d-m_1}=0,
$$
and be a 2-term one.
\begin{flushright}
\bf QED
\end{flushright}\rm

\subsection{The exceptional birational equivalences}
\label{lowdeg}
\subsubsection{Cubics}
\label{cubics}
For $d=3$ we have
\begin{center}
\begin{tabular}{|c|c|c|}
\hline
Name& Equation & j-invariant \\
\hline
Fermat (diagonal) type&  $x^3+y^3+z^3=0$ & 0  \\
Small Jordan type& $xy^2+y^3+z^3=0$ & 0  \\
block type&  $x^3+yz^2+y^2z=0$ & 0 \\
Big Jordan type& $x^3+xy^2+yz^2=0$& 1728 \\
Klein type &  $xy^2+yz^2+zx^2=0$  & 0 \\
 \hline
\end{tabular}             \end{center}
so we have only two birational types of 3-term curves.

\
\subsubsection{Quartics}
\label{quartics}
For $d=4$ we note that the curve $x^4+y^4+z^4=0$ is projectively equivalent to $u^4+wv^3+vw^3=0$.
$$
\left( \begin{array}{ccc}
1& 0 & 0 \\
0 & \frac{2^{\frac{1}{4}}}{2} & \frac{2^{\frac{3}{4}}}{4}+i \frac{2^{\frac{3}{4}}}{4} \\
0 & \frac{2^{\frac{1}{4}}}{2} & \frac{2^{\frac{3}{4}}}{4}+i \frac{2^{\frac{3}{4}}}{4} \end{array} \right)
\left( \begin{array}{c}
u \\
v  \\
w  \end{array} \right)
=
\left( \begin{array}{c}
x \\
y  \\
z  \end{array} \right),
$$
so there are four birational types of 3-term curves.

\section{Diagonal automorphisms}
\subsection{Determinants}
Specify the power matrices and calculate the determinants.
\begin{center}
\begin{tabular}{|c|c|c|c|}
\hline
Name& Equation &  Matrix & Determinant\\
\hline Fermat  &  $x^d+y^d+z^d=0$ & $\left (
 \begin{array}{ccc}
  d&0&0\\
  0&d&0\\
  0&0&d\\
 \end{array}
\right )$&$d^3$  \\

small Jordan & $x^d+y^d+yz^{d-1}=0$ &  $\left (
 \begin{array}{ccc}
  d&0&0\\
  0&d&0\\
  0&1&d-1\\
 \end{array}
\right )$ &$d^2(d-1)$\\
block &  $x^d+y^{d-1}z+yz^{d-1}=0$ & $\left (
 \begin{array}{ccc}
  d&0&0\\
  0&d-1&1\\
  0&1&d-1\\
 \end{array}
\right )$ &$d^2(d-2)$\\
big Jordan & $x^d+xy^{d-1}+yz^{d-1}=0$ & $ \left (
 \begin{array}{ccc}
  d&0&0\\
  1&d-1&0\\
  0&1&d-1\\
 \end{array}
\right )$&$d(d-1)^2$ \\
Klein &  $xy^{d-1}+yz^{d-1}+zx^{d-1}=0$
  &  $\left (
 \begin{array}{ccc}
  1&d-1&0\\
  0&1&d-1\\
  d-1&0&1\\
 \end{array}
\right )$&$d(d^2-3d+3)$\\
 \hline
\end{tabular}
\end{center}
\subsection{The  group of diagonal automorphisms}
It consists, by definition, of the projective linear automorphisms
of the form
$$
(x:y:z)\mapsto(\rho^Ax:\rho^By:\rho^Cz),
$$
where $A,B,C$ are defined up to a common summand. Obviously the
triple $A,B,C$ defines an automorphism of the corresponding 3-term
curve if \eq{ P \left (\begin{array}{c} A\\B\\C\\ \end{array} \right
) = \left( \begin{array}{c}
0 \\
0  \\
0  \end{array} \right). } Since all the five determinants divide
$d$, considering the last equation over
$\frac{\mathbb{Z}}{d\mathbb{Z}}$, we see that it has nonzero
solutions, hence all our curves have non-trivial diagonal
automorphisms. Moreover, the order of this group is
$\frac{\Delta}{d}$.

\section{Belyi functions}
\label{Belyi}
For all five types of curves the Belyi functions are given by the
following table
%\eql{table2}{

\begin{center}
\begin{tabular}{|c|c|c|}
\hline
curve &  Belyi \ function  & the \ degree \ of \ Belyi \ function  \\
\hline $x^d+y^d+z^d=0$                           & $\br{\dfrac{x}{z}}^d$           & $d^2$ \\
$x^d+y^d+yz^{d-1}=0$                   &  $\br{\dfrac{y}{z}}^{d-1}+1$ & $d(d-1)$ \\
$x^d+zy^{d-1}+yz^{d-1}=0$          &  $\br{\dfrac{y}{z}}^{d-2}+1$ & $d(d-2)$ \\
$x^d+zy^{d-1}+xz^{d-1}=0$          & $\br{\dfrac{x}{z}}^{d-1}+1$ & $(d-1)^2$ \\
$x^{d-1}y+y^{d-1}z+z^{d-1}x=0$ & $\dfrac{x^{d-1}}{zy^{d-2}}$ & $d^2-d$\\
 \hline
\end{tabular}
\end{center}
So all our curves have remarkably low Belyi heights.

\section*{Acknowledgements}
The authors are grateful to Natalia Amburg, Petr Dunin-Barkowski, Andrey Smirnov and Ilya Zhdanovsky for fruitful discussions and very helpful remarks.  This work is partly supported by the Russian  President's Grant of Support for the Scientific Schools NSh-3035.2008.2, by joint grant 09-01-92440-CE, by RFBR grant 07-01-00441 (G.B. Shabat) and by RFBR grant 07-02-00878 (A. Sleptsov).

\end{document}